\documentclass[preprint,12pt]{elsarticle}

\usepackage{amssymb}
\usepackage{amsmath}
\usepackage{color}
\usepackage{amsthm}

\journal{Information Sciences}

\begin{document}

\begin{frontmatter}



\title{Conditional investment strategy in evolutionary trust games with repeated group interactions}


\author[label1,label2]{Linjie Liu}
\author[label2]{Xiaojie Chen\corref{cor}}
\cortext[cor]{Corresponding author} \ead{xiaojiechen@uestc.edu.cn}
\address[label1]{College of Science, Northwest A \& F University, Yangling 712100, China}
\address[label2]{School of Mathematical Sciences, University of
Electronic Science and Technology of China, Chengdu 611731, China}

\begin{abstract}
It has a long tradition to study trust behavior among humans or artificial agents by investigating the trust game. Although previous studies based on evolutionary game theory have revealed that trust and trustworthiness can be promoted if network structure or reputation is considered, they often assume that interactions among agents are one-shot and investors do not consider the investment environment before making decisions, which collide with many realistic situations. In this paper, we introduce the conditional investment strategy into the repeated $N$-player trust game, in which conditional investors decide to invest or not depending on their assessment of the trustworthiness level of the group. By using the approach of the Markov decision process, we study the evolutionary dynamics of trust in repeated group interactions with the conditional investment strategy. We find that conditional investors can form an effective alliance with trustworthy trustees, hence they can sweep out untrustworthy trustees. Moreover, we verify that such alliance can exist in a wide range of model parameters. These results may explain why trusting in others and reciprocating them with trustworthy actions can be sustained in game interactions among intelligent agents.
\end{abstract}


\end{frontmatter}


\section{Introduction}

\noindent Trust is a value-laden concept in Distributed Artificial Intelligence (DAI), and it is also fundamental for cooperative behavior and other forms of prosocial behavior \cite{Yu2015,liukbs20,Perc2017,azolnki2014interface,Urena2019is}. 
However, a basic fact accepted by the public is that trust is often accompanied by risks. 
As Deutsch defined, entering a trust relationship is to choose an ambiguous path that can lead to a beneficial event or a harmful event depending on the behavior of the other person \cite{deutsch1962cooperation}. 
In other words, there are risks in the interactions with potentially untrustworthy agents, while the interactions with trusted agents will lead to benefits \cite{swinth1967establishment}. 
Thus understanding how trust and trustworthiness evolve in Multi-agent Systems (MAS) and exploring the conditions in which they will emerge have always been great challenges \cite{han2021csr}.

Evolutionary game theory provides a theoretical framework to study the above problems~\citep{fang2020is,Liangis2021,Johnson2015,perc13inter,szonoki2013PRX,li2014}, 
and the trust game (TG) has been widely used to study trust and trustworthiness as a typical paradigm \citep{abbass2015n,hu2021adaptive}. 
The classical TG model involves interactions between an investor and a trustee where the investor first decides whether he/she is willing to invest his/her funds to the trustee, 
and then the trustee decides how much to give back to the investor (or not at all) \citep{berg1995trust,king2005getting}. 
Previous theoretical work has proved that the subgame perfect equilibrium of TG is that the investor invests zero and the trustee returns zero \citep{guth1997cooperation}. 
The above two-player TG model perfectly characterizes pairwise interactions between an investor and a trustee and has recently been extended to group interactions of multiple players \cite{abbass2015n,chica2017networked,hu2021adaptive,kumar2020evolution}. 
A representative $N$-player trust game (NTG) framework was developed by Abbass \emph{et al.} \cite{abbass2015n}, and then introduced to a population of agents playing TG concurrently in a well-mixed environment. 
Along this line, many researchers have considered network structure \citep{chica2017networked,chica2019effects,kumar2020evolution}, punishment strategy \citep{fang2021epjb}, asymmetric
demographic parameters \cite{lim2020stochastic}, and reputation \citep{hu2021adaptive}, into the NTG model to explore the evolution of trust.

Although there were attempts for enhancing the level of trust from different aspects, we notice that these studies are usually carried out in the framework of one-shot interaction, to our knowledge \cite{chica2017networked,chica2019effects,fang2021epjb,kumar2020evolution,lim2020stochastic}. Indeed, real interactions are repeated rather than the one-shot interaction assumed above \cite{Chong05,Chiong12,martinez2020signalling,quan_21_kbs,liu2022jrsi}. For example, we are likely to interact frequently with friends, co-workers, and economic partners in our daily life. Second, agents' behaviors may be affected by noise (such as imitation error and behavioral mutation) in the process of social learning \cite{Sigmund2010,sun2021combination,vasconcelos2013bottom}, which is usually ignored by previous theoretical models.
Furthermore, in repeated interactions, investors have the opportunity to adjust their actions according to the investment environment they encountered \cite{martinez2020signalling,quan_19_chaos,liu2022jrsi}. Furthermore, it is very common for investors to choose conditional investment behavior in their daily economic activities when interacting with other agents \cite{Ahmed2011,Bond2013}. Concretely, the decision to invest or not will depend on one's own assessment of the abundance of trustworthy trustees in the group. However, it is still unclear what the effects of conditional investment strategy on the evolution of trust in the repeated NTG are.

\begin{figure}[!h]
\begin{center}
\includegraphics[width=4in]{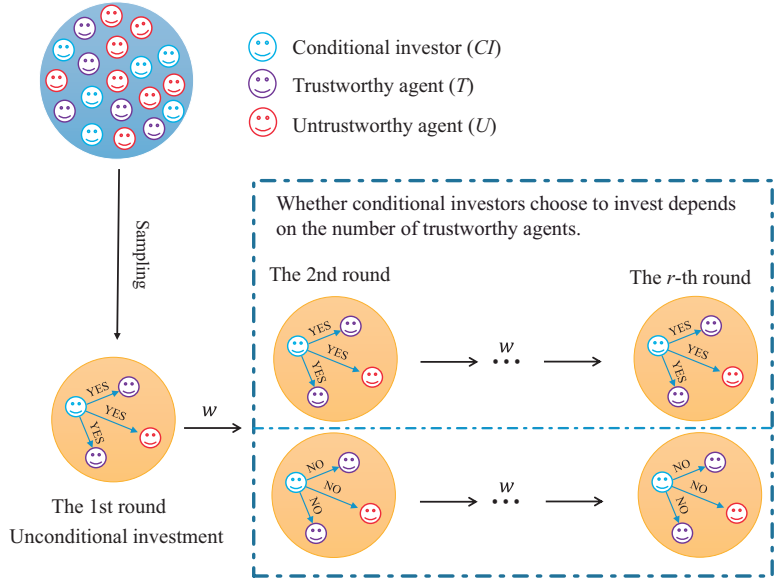}
\caption{Illustration of two different interactions scenarios in repeated $N$-player trust games where the repeated probability is set as $w$. $N$ agents are randomly sampled to form a group for playing the repeated trust game where a conditional investor will act as an unconditional investor in the first round, and then can agree or refuse to invest depending on whether the number of trustworthy agents in the group exceeds the expected threshold in subsequent game rounds.}\label{fig1}
\end{center}
\end{figure}

In order to answer the above question, we construct an evolutionary game model based on the NTG where the interactions among agents are repeated (see Fig.\ref{fig1}). In our model, we consider three strategists: investor, untrustworthy trustee, and trustworthy trustee. Different from previous works \cite{hu2021adaptive,kumar2020evolution}, we assume investors will invest deterministically in the first round of the game, but in the subsequent game rounds their decision to invest or not depends on the number of the trustworthy trustee in the group. Using the approach of the Markov decision process (MDP), we study the stochastic evolutionary dynamics in finite populations. We find that the introduction of the conditional investment strategy can lead to the emergence of high trust and high trustworthiness in repeated group interactions. Furthermore, we verify that such evolutionary outcomes are robust against model parameters.

\section{Related Works}

The formalization of the NTG that we follow here was proposed in Ref.~\cite{abbass2015n}, where Abbass \emph{et al.} investigated the evolutionary dynamics of trust in an infinite well-mixed population. By analyzing the replicator dynamics, they found that both the whole society and all agents can obtain the maximal wealth when the initial population contains no untrustworthy agents. While the initial population consists of one single untrustworthy agent, untrustworthy agents would spread rapidly to the whole population. Furthermore, they found that the system will eventually converge to a stable state, in which untrustworthy agents will not completely occupy the whole population, because a fraction of the population would always remain trustworthy even if there are few or no investors.

Nevertheless, the analysis of replicator dynamics in the previous literature requires that the network structure is a complete graph, and thus this approach cannot be used to analyze the evolutionary dynamics when agents interact with local neighborhoods in social networks. Subsequently, Chica \emph{et al.} studied the evolution of trust in different social networks including regular lattices, scale-free, and random networks \cite{chica2017networked}. By performing simulations, they found that trust can be promoted when agents interact on a social network even if there are untrustworthy agents in the population initially. In addition, the level of trust is influenced by network structure, the temptation to defect, and the initial number of untrustworthy trustees in the population. Along this line, the effects of different evolutionary update rules on promoting the evolution of trust were investigated \cite{chica2019effects}, and simulation results showed that updating rules play an important role in promoting trust and improving global net wealth.

However, some other studies have found that social networks cannot produce a high level of trust \cite{hu2021adaptive,kumar2020evolution}. Recently, Kumar \emph{et al.} investigated the trust game where the trustor's investment and the trustee's return of the investment are two important parameters on different social networks \cite{kumar2020evolution}. By performing Monte Carlo simulations, they found that the network structure has little effect on the evolution of
trust and trustworthiness. In particular, trust cannot evolve in well-mixed populations, lattices, random, or scale-free networks. Considering that in the real society, the investment behavior of investors is often related to the reputation of the trustees, and agents with good reputation are more likely to attract investors' investment, Hu \emph{et al.} studied networked NTG with an adaptive reputation based on the third-party moral assessment system \cite{hu2021adaptive}. They showed that the frequency of untrustworthy trustees will decrease when rational investors can get the reputation score of the trustee.

The literature mentioned above does not consider the case in which agents' decision-making may be affected by random factors, including imitation error and behavioral mutation. Considering this, Lim explored the evolutionary dynamics of trust in finite populations and found that the combination of strong selection in the population of investors and weak selection in the population of trustees can promote the emergence of high trustworthiness and high trust \cite{lim2020stochastic}. It is worth noting that they analyzed the stationary distribution of the system in the assumption of weak-mutation limit. The above assumption has been justified in population genetics \cite{crow1970,Ewens2004}, but this approximate does not seem to be suitable for modeling social learning \cite{vasconcelos2013bottom}.

Different from previous works both in infinite populations \cite{abbass2015n} and structure populations \cite{chica2017networked}, we study the evolutionary dynamics of trust in a finite population where imitation error and behavioral mutation are both exist. Besides, we release the assumptions that the mutation rate is significantly small \cite{lim2020stochastic} and the interactions between agents are one-shot \cite{abbass2015n,chica2017networked,kumar2020evolution,lim2020stochastic}, and investigate the stochastic
dynamics of trust in finite populations with repeated group interactions when mutation rates are arbitrarily large. Accordingly, we are committed to solving the following questions: can potentially richer evolutionary dynamics be produced? Can a high level of trust and trustworthiness be reached in the NTG with repeated group interactions? To answer these questions, we establish an NTG model where the conditional investment strategy is involved in repeated group interactions. Concretely, conditional investors decide to invest only when the number of trustworthy trustees in the game group reaches their expected threshold. We study the evolutionary dynamics of strategies in a finite population by using the approach of MDP. Our concrete theoretical model and methods are presented in the following section.

\section{Theoretical Model and Methods}
\subsection{Repeated NTG}

Let us assume a finite well-mixed population with $Z$ agents who play a repeated NTG \cite{abbass2015n,chica2017networked}.
At every time step, $N$ agents are chosen and offered the opportunity to participate in the repeated NTG. As considered in previous work \cite{liu2022jrsi}, the trust game will be repeated with probability $w$ with $0<w<1$ (also called the discount factor in previous work \cite{Hilbe2018}), resulting in an average number of $r=1/w$ rounds. In this game, each agent needs to make two decisions in advance. First, agent needs to decide whether to act as a trustee or an investor. The second is to decide whether to be trustworthy or not. Here, we consider three baseline strategists in this NTG:

{\begin{itemize}
\item Conditional Investor ($CI$): Pays an observation cost $\sigma$ before the game to gather information including the number of trustworthy agents and untrustworthy agents in the group, and always chooses to invest $tv$ in the first round. Subsequently, $CI$ agents invest only when the number of trustworthy agents in the group is not less than their expected threshold $M$.
\item Trustworthy Agent ($T$): Returns the received fund obtained in each round multiplied by $R_{T}$ to the investors.
\item Untrustworthy Agent ($U$): Returns nothing to the investors in all game rounds.
\end{itemize}

We shall first notice that, by adopting $CI$, an agent will always invest in the first round and subsequently either act as an investor or as an outsider, depending on whether the number of trustworthy trustees in the group has reached the expected level $M$ (see Fig.\ref{fig1}). If the expected level is reached, $CI$ agents are willing to pay $tv$ to trustees, and each trustee receives the same amount of benefit $\frac{N_{CI}tv}{N_{T}+N_{U}}$, where $N_{CI}, N_{T}$, and $N_{U}$ respectively denote the number of $CI, T$, and $U$ agents in the group. Then each $T$ agent returns $\frac{R_{T}tv}{N_{T}+N_{U}}$ to every $CI$ agent and reserves $\frac{R_{T}N_{CI}tv}{N_{T}+N_{U}}$ for himself/herself. While each $U$ agent returns nothing to $CI$ agents, but keeps $\frac{R_{U}N_{CI}tv}{N_{T}+N_{U}}$ for himself/herself. The parameters $R_{T}$ and $R_{U}$ respectively denote the multiply factor of $T$ agents and $U$ agents.
If the number of $T$ agents in the group does not reach the expected level, every agent in the group can obtain nothing from the subsequent game rounds. Thus the payoffs of $CI, T$, and $U$ agents obtained from the game are depicted by

\begin{align}
\Pi_{CI}&=\left\{
\begin{aligned}\label{eq1}
&\frac{R_{T}N_{T}tv}{N-N_{CI}-1}-tv+(\frac{R_{T}N_{T}tv}{N-N_{CI}-1}-tv)(r-1)\Theta(N_{T}-M)-\sigma,  \quad  \text{if}\ N_{CI}\neq N-1\ ;\\
&0,     \quad   \text{otherwise.}
\end{aligned}
\right.\\
\Pi_{T}&=\left\{
\begin{aligned}
&\frac{R_{T}N_{CI}tv}{N-N_{CI}}+\frac{R_{T}N_{CI}tv}{N-N_{CI}}(r-1)\Theta(N_{T}+1-M), \quad  \text{if}\ N_{CI}\neq N\ ;\\
&0,    \quad    \text{otherwise.}
\end{aligned}
\right.\\
\Pi_{U}&=\left\{
\begin{aligned}\label{eq3}
&\frac{R_{U}N_{CI}tv}{N-N_{CI}}+\frac{R_{U}N_{CI}tv}{N-N_{CI}}(r-1)\Theta(N_{T}-M), \quad  \text{if}\ N_{CI}\neq N\ ;\\
&0,    \quad    \text{otherwise,}
\end{aligned}
\right.
\end{align}
\noindent where $\Theta(k)$ is the Heaviside function (that is, $\Theta(k)=1$ whenever
$k\geq0$, being zero otherwise). The threshold values $M$ range from 1 to $N-1$. If $M<1$ all $CI$ agents will choose to invest unconditionally and $M>N-1$ means always choosing not to invest fund.

In a finite population with $Z$ agents, $N$ agents are randomly sampled to form a game group and accumulate their payoffs by interacting with other group agents according with the payoff functions of NTG defined above. Concretely, when there are $i_{CI}$ agents choosing to adopt the $CI$ strategy, $i_{T}$ agents choosing to adopt the $T$ strategy, $i_{U}=Z-i_{CI}-i_{T}$ agents choosing to adopt the $U$ strategy in a finite population, the probability of finding $j_{CI}$ $CI$ agents, $j_{T}$ $T$ agents, and $j_{U}$ $U$ agents in a sample game group can be depicted by the multiple hypergeometric distribution, given as
\begin{eqnarray*}
H(j_{CI},j_{T},N,i_{CI},i_{T},Z)=\frac{\tbinom{i_{CI}}{j_{CI}}\tbinom{i_{T}}{j_{T}}\tbinom{Z-i_{CI}-i_{T}}{N-j_{CI}-j_{T}}}{\tbinom{Z}{N}}, 
\end{eqnarray*}
which describes the configuration of the population for a given time. Accordingly, for a given configuration $\textbf{i}=({i_{CI},i_{T}})$, the average payoffs of $CI, T$, and $U$ agents can be respectively computed as
\begin{eqnarray*}
f_{CI}&=&\sum\limits_{j_{CI}=0}^{N-1}\sum\limits_{j_{T}=0}^{N-1-j_{CI}}\frac{\tbinom{i_{CI}-1}{j_{CI}}\tbinom{i_{T}}{j_{T}}\tbinom{Z-i_{CI}-i_{T}}{N-1-j_{CI}-j_{T}}}{\tbinom{Z-1}{N-1}}\\
&\times&\Pi_{CI}(j_{CI}+1,j_{T}),\\
f_{T}&=&\sum\limits_{j_{CI}=0}^{N-1}\sum\limits_{j_{T}=0}^{N-1-j_{CI}}\frac{\tbinom{i_{CI}}{j_{CI}}\tbinom{i_{T}-1}{j_{T}}\tbinom{Z-i_{CI}-i_{T}}{N-1-j_{CI}-j_{T}}}{\tbinom{Z-1}{N-1}}\\
&\times&\Pi_{T}(j_{CI},j_{T}+1),\\
f_{U}&=&\sum\limits_{j_{CI}=0}^{N-1}\sum\limits_{j_{T}=0}^{N-1-j_{CI}}\frac{\tbinom{i_{CI}}{j_{CI}}\tbinom{i_{T}}{j_{T}}\tbinom{Z-i_{CI}-i_{T}-1}{N-1-j_{CI}-j_{T}}}{\tbinom{Z-1}{N-1}}\\
&\times&\Pi_{U}(j_{CI},j_{T}),
\end{eqnarray*}
where $\Pi_{CI}, \Pi_{T},$ and $\Pi_{U}$ are respectively shown in Eqs. (\ref{eq1})-(\ref{eq3}).

\subsection{Social Learning}

The average payoffs of agents obtained above stand for their social success, then we can analyze the evolutionary dynamics of strategy adopting by using social learning \cite{quan_17_SR,quan_19_SR,Sigmund2010}, which means that the most successful strategy will more often tend to be adopted by other agents. Generally, the above social learning can be characterized by the pairwise comparison rule. Concretely, the probability that an agent adopting strategy $A$ with the payoff $f_{A}$ imitates another agent's strategy $B$ with payoff $f_{B}$ is given by the Fermi function~\cite{szabo1998evolutionary}
\begin{eqnarray}
P(f_{B}-f_{A})=\frac{1}{1+\exp(-\beta(f_{B}-f_{A}))},
\end{eqnarray}
where $\beta$ characterizes the intensity of selection, controlling how the imitation process depends on the difference between the payoffs of two agents. For $\beta \rightarrow \infty$ (strong imitation), any difference in the payoffs will produce a sizeable effect on imitation probability. For $\beta \rightarrow 0$ (weak selection), the strategy is imitated randomly. In between these extremes, the difference of the payoffs and stochastic effects associated with errors can both impact imitation probability. 

Furthermore, we introduce the exploration term: with probability $1-\mu$, an agent with strategy $A$ adopts another agent's strategy $B$ according to the social learning rule described above, and with probability $\mu$, he/she randomly selects a strategy different from the current one from the strategy space. Thus, the probability that an agent with strategy $A$ adopts strategy $B$ according to the mutation-selection process can be written as
\begin{eqnarray}\label{transeq}
T_{A\rightarrow B}&=&(1-\mu)[\frac{i_{A}}{Z}\frac{i_{B}}{Z-1}\frac{1}{1+\exp({-\beta(f_{B}-f_{A})})}]+\mu\frac{i_{A}}{2Z}.
\end{eqnarray}
The existence of behavioral mutation makes the population never fixate in any of the three possible monomorphic configurations. This fact renders the stationary distribution and the gradient of selection as the appropriate quantities to analyze the behavior of the population.

\subsection{Markov Decision Process}

We know that the update process only relies on the current state $\textbf{i}(t)$ of the system, so that $\textbf{i}(t)=\{i_{CI}, i_{T}\}$ has Markov property. Thus the evolutionary dynamics of \emph{CI}, \emph{T}, and \emph{U} can be described by the Markov chain in two-dimensional space. Then the evolutionary dynamics of the system can be analyzed by investigating the probability distribution function $p_{\textbf{i}}(t)$ providing information about the pervasiveness of each configuration at time \emph{t}, and it satisfies the following discrete time Master Equation \cite{van1992stochastic}
\begin{eqnarray*}
p_{\textbf{i}}(t+\tau)-p_{\textbf{i}}(t)=\sum_{\textbf{i}^{'}}\left\{T_{\textbf{i}\textbf{i}^{'}}p_{\textbf{i}^{'}}(t)-T_{\textbf{i}^{'}\textbf{i}}p_{\textbf{i}}(t)\right\},&
\end{eqnarray*}
where $T_{\textbf{i}\textbf{i}^{'}}$ and $T_{\textbf{i}^{'}\textbf{i}}$ denote the transition probabilities between configurations $\textbf{i}'$ and $\textbf{i}$. Technically, we can obtain the so-called stationary distribution $\bar{p}_{\textbf{i}}$, by searching the eigenvector associated with the eigenvalue 1 of the transition matrix $\Lambda=[T_{\textbf{ij}}]^{T}$ with dimension $\frac{(Z+1)(Z+2)}{2}\times \frac{(Z+1)(Z+2)}{2}$. The transition probability between two adjacent states can be computed as
\begin{eqnarray*}
T_{\textbf{i}(i_{CI}, i_{T})\rightarrow \textbf{i}^{'}(i_{CI}+1, i_{T})}&=&T_{U\rightarrow CI},\\
T_{\textbf{i}(i_{CI}, i_{T})\rightarrow \textbf{i}^{'}(i_{CI}-1, i_{T})}&=&T_{CI\rightarrow U},\\
T_{\textbf{i}(i_{CI}, i_{T})\rightarrow \textbf{i}^{'}(i_{CI}, i_{T}+1)}&=&T_{U\rightarrow T},\\
T_{\textbf{i}(i_{CI}, i_{T})\rightarrow \textbf{i}^{'}(i_{CI}, i_{T}-1)}&=&T_{T\rightarrow U},\\
T_{\textbf{i}(i_{CI}, i_{T})\rightarrow \textbf{i}^{'}(i_{CI}-1, i_{T}+1)}&=&T_{CI\rightarrow T},\\
T_{\textbf{i}(i_{CI}, i_{T})\rightarrow \textbf{i}^{'}(i_{CI}+1, i_{T}-1)}&=&T_{T\rightarrow CI}.
\end{eqnarray*}
The transition probability between two nonadjacent states $\textbf{i}$ and $\textbf{v}$ is $T_{\textbf{i} \textbf{v}}=0$. Thus the probability of the system staying in the current state is
\begin{eqnarray*}
T_{\textbf{i}(i_{CI}, i_{T})\rightarrow \textbf{i}(i_{CI}, i_{T})}&=&1-\sum_{\textbf{i}\neq\textbf{i}^{'}}(T_{\textbf{i}\textbf{i}^{'}}+T_{\textbf{i}^{'}\textbf{i}}).
\end{eqnarray*}
In one discrete time step, the probability of the system transferring from one state to another can be calculated by equation (\ref{transeq}).

In addition to the analysis of the stationary distribution of the system, another important quantity for studying the evolutionary dynamics in finite populations
is the gradient of selection, which indicates the most likely evolutionary path when the system leaves the current configuration \cite{vasconcelos2013bottom}.
The gradient of selection is described as
\begin{eqnarray}\label{grad}
\nabla_{\textbf{i}}&=(T_{\textbf{i}}^{CI+}-T_{\textbf{i}}^{CI-})\textbf{u}_{\textbf{CI}}+(T_{\textbf{i}}^{T+}-T_{\textbf{i}}^{T-})\textbf{u}_{\textbf{T}},&
\end{eqnarray}
where $\textbf{u}_{\textbf{CI}}$ and $\textbf{u}_{\textbf{T}}$ are a set of standard orthogonal bases, and we set $\textbf{u}_{\textbf{CI}}=(1,0)^{T}$ and $\textbf{u}_{\textbf{T}}=(0,1)^{T}$ in this work. $T_{\textbf{i}}^{CI+} (T_{\textbf{i}}^{CI-})$ and $T_{\textbf{i}}^{T+} (T_{\textbf{i}}^{T-})$ respectively denote the probabilities that the numbers of $CI$ agents and $T$ agents increase (decrease) one, which read
\begin{eqnarray*}
T_{\textbf{i}}^{CI+}&=& T_{U\rightarrow CI}+T_{T\rightarrow CI},\\
T_{\textbf{i}}^{CI-}&=&T_{CI\rightarrow U}+T_{CI\rightarrow T},\\
T_{\textbf{i}}^{T+}&=& T_{U\rightarrow T}+T_{CI\rightarrow T},\\
T_{\textbf{i}}^{T-}&=&T_{T\rightarrow U}+T_{T\rightarrow CI}.
\end{eqnarray*}
Furthermore, we provide an important index to describe the average level of each strategy. Concretely, the average levels of $CI, T,$ and $U$ strategies, averaging over all possible states $\textbf{i}$, weighted with the corresponding stationary distribution $\bar{p_{\textbf{i}}}$, are  computed as
\begin{eqnarray*}
\bar{\rho}_{CI}&=&\sum_{\textbf{i}}\frac{\textbf{i}_{i_{CI}}\bar{p_{\textbf{i}}}}{Z},\\
\bar{\rho}_{T}&=&\sum_{\textbf{i}}\frac{\textbf{i}_{i_{T}}\bar{p_{\textbf{i}}}}{Z},\\
\bar{\rho}_{U}&=&\sum_{\textbf{i}}\frac{\textbf{i}_{i_{U}}\bar{p_{\textbf{i}}}}{Z},
\end{eqnarray*}
where $\textbf{i}_{i_{CI}}, \textbf{i}_{i_{T}}$, and $\textbf{i}_{i_{U}}$ denote the number of $CI, T$, and $U$ agents in the configuration $\textbf{i}$, respectively.

Subsequently, by using $f_{S}(\textbf{i})$ combined with the stationary distribution $\bar{p}_{\textbf{i}}$, we can calculate the average payoff of one agent with $S$ strategy as
\begin{eqnarray*}
\bar{f}_{S}=\sum_{\textbf{i}}\bar{p_{\textbf{i}}}f_{S}(\textbf{i}),
\end{eqnarray*}
where $S = CI, T,$ or $U$.

In the following, we investigate the gradient of selection and the stationary distribution to study the evolutionary dynamics
of $CI$, $T$, and $U$ strategies in finite well-mixed populations.

\begin{figure*}[!t]
\begin{center}
\includegraphics[width=5.5in]{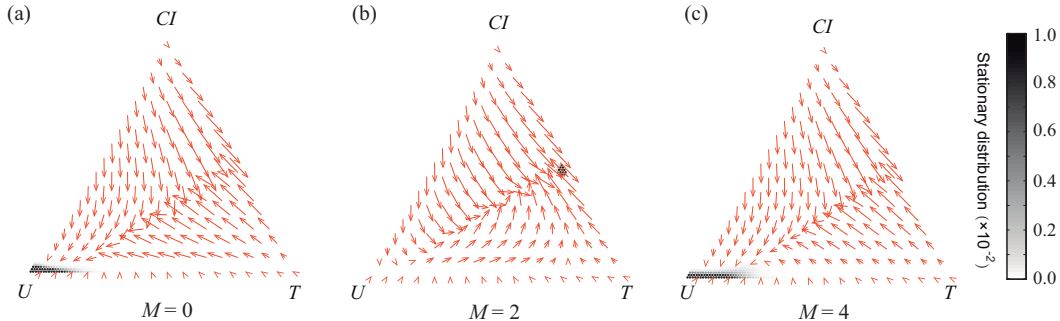}
\caption{Evolutionary dynamics of Conditional Investor ($CI$), Trustworthy agents ($T$), and Untrustworthy agents ($U$) for three different threshold values $M$ in finite populations. Arrows in the simplex $S_{3}$ represent the gradient of selection, which shows the most likely evolutionary paths after the system leaving the current state, calculated from equation (\ref{grad}). The colorbar describes the values of stationary distribution at each configuration. The darker the color, the longer the population spends in these states. Parameters are $Z=100, N=4, tv=1, R_{U}=8, R_{T}=6, \mu=1/Z, w=0.8, \sigma=0.1,$ and $\beta=5$. Here, $M=0$ in panel (a), $M=2$ in panel (b), and $M=4$ in panel (c).}\label{fig2}
\end{center}
\end{figure*}

\begin{figure}[!t]
\begin{center}
\includegraphics[width=4in]{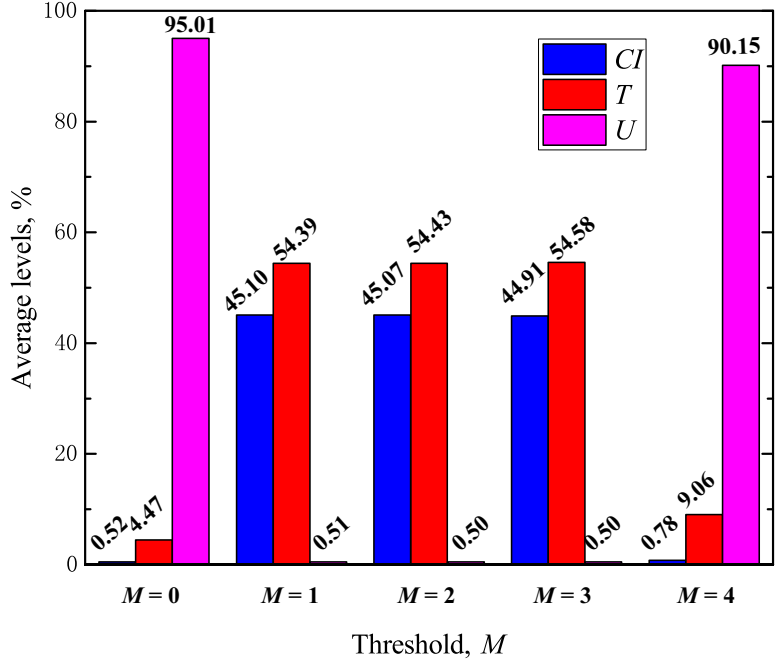}
\caption{Average frequency of strategies $CI$ (blue cylinder), $T$ (red cylinder), and $U$ (purple cylinder) for different investment threshold $M$. The intermediate investment threshold is conducive to the formation of $CI$ and $T$ alliance, and then effectively improves the level of trust. Parameters are $Z=100, N=4, tv=1, R_{U}=8, R_{T}=6, w=0.8, \mu=1/Z, \sigma=0.1,$ and $\beta=5$.}\label{fig3}
\end{center}
\end{figure}

\begin{figure}[!t]
\begin{center}
\includegraphics[width=4in]{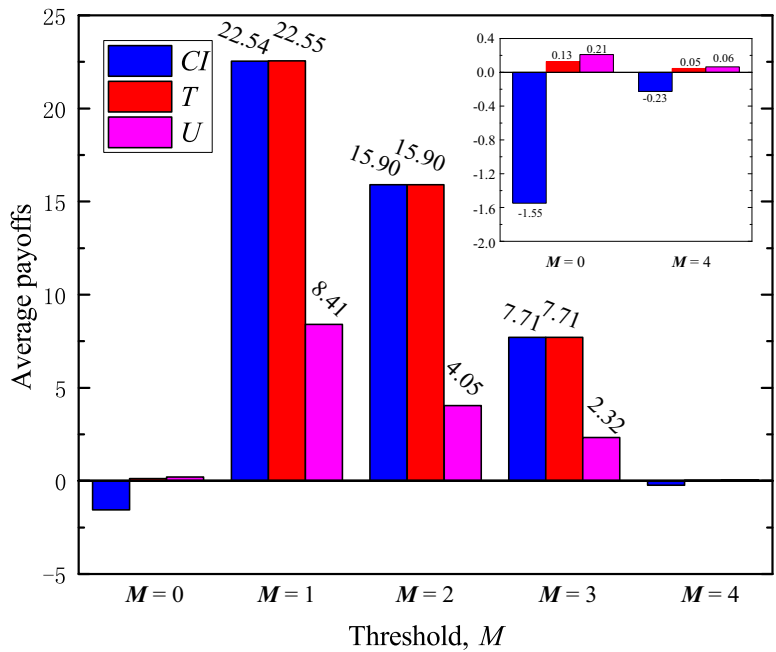}
\caption{Average payoffs of $CI$ (blue cylinder), $T$ (red cylinder), and $U$ (purple cylinder) agents for different investment threshold values $M$. The advantage of the alliance formed by $CI$ and $T$ agents can be boosted by the intermediate investment threshold $M$. Parameters are $Z=100, N=4, tv=1, R_{U}=8, R_{T}=6, \mu=1/Z, w=0.8, \sigma=0.1,$ and $\beta=5$.}\label{fig4}
\end{center}
\end{figure}

\section{Results}

\begin{figure*}[!h]
\begin{center}
\includegraphics[width=5.5in]{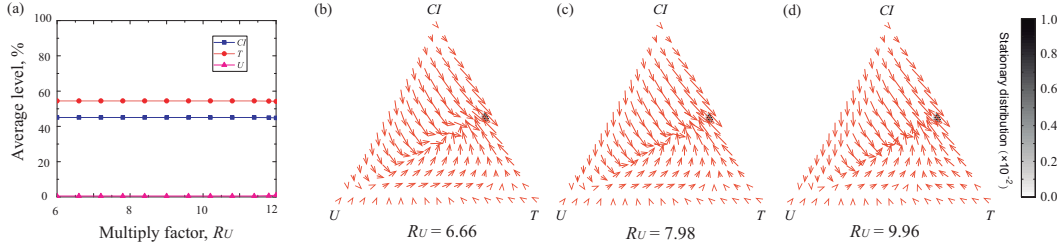}
\caption{Robustness of the evolutionary outcomes of $CI, T$, and $U$ to the changes of the multiply factor of $U$ agents $R_{U}$. Panel (a) shows the average levels of three strategies as a function of $R_{U}$. Panels (b)-(d) show the stationary distribution and the gradient of selection for three different $R_{U}$ values. Parameter values are $Z=100, N=4, M=2, tv=1, w=0.8, R_{T}=6, \mu=1/Z, \beta=5, \sigma=0.1$ in panel (a); $Z=100, N=4, M=2, tv=1, R_{U}=6.66, w=0.8, R_{T}=6, \mu=1/Z, \beta=5, \sigma=0.1$ in panel (b); $Z=100, N=4, M=2, tv=1, R_{U}=7.98, w=0.8, R_{T}=6, \mu=1/Z, \beta=5, \sigma=0.1$ in panel (c); $Z=100, N=4, M=2, tv=1, R_{U}=9.96, w=0.8, \sigma=0.1, R_{T}=6, \mu=1/Z, \beta=5$ in panel (d).}\label{fig5}
\end{center}
\end{figure*}

\begin{figure*}[!h]
\begin{center}
\includegraphics[width=5.5in]{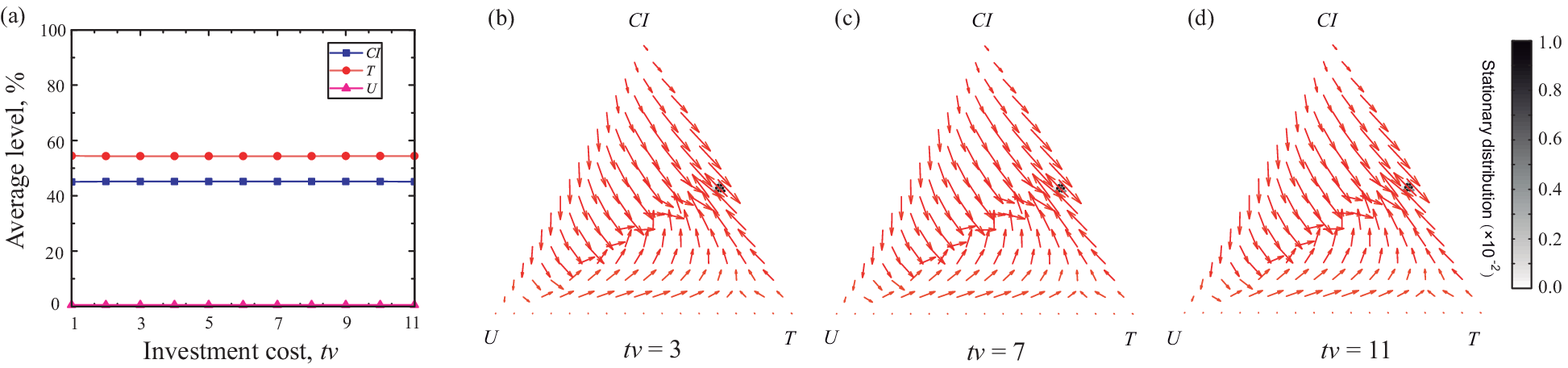}
\caption{Robustness of the evolutionary outcomes of $CI, T$, and $U$ to the changes of the investment cost $tv$. Panel (a) shows the average levels of three strategies as a function of $tv$. Panels (b)-(d) show the stationary distribution and the gradient of selection for three different $tv$ values. Parameter values are $Z=100, N=4, M=2, R_{U}=8, w=0.8, R_{T}=6, \mu=1/Z, \beta=5, \sigma=0.1$ in panel (a); $Z=100, N=4, M=2, tv=3, R_{U}=8, w=0.8, R_{T}=6, \mu=1/Z, \beta=5, \sigma=0.1$ in panel (b); $Z=100, N=4, M=2, tv=7, R_{U}=8, w=0.8, R_{T}=6, \mu=1/Z, \beta=5, \sigma=0.1$ in panel (c); $Z=100, N=4, M=2, tv=11, R_{U}=8, w=0.8, \sigma=0.1, R_{T}=6, \mu=1/Z, \beta=5$ in panel (d).}\label{fig6}
\end{center}
\end{figure*}

We first present the results of evolutionary dynamics in finite populations with the conditional investment strategy for
different values of $M$. In Fig. \ref{fig2}, we investigate the stationary distribution and the gradient of selection to study what roles conditional investment plays in the evolution of trust. Clearly, if the investment tendency is strong enough ($M=0$), $CI$ agents will become unconditional investors. We can notice that the population will spend a significant time on configurations near the $TU$-edge ($T$ means trustworthy agents and $U$ means untrustworthy agents) of the triangle simplex, as shown in Fig. \ref{fig2}(a). Besides, most of the arrows flow to the coexistence states of $T$ and $U$ where $U$ is highly prevalent. The weakening of investment tendency will lead the population to spending most of the time in states with a high prevalence of $T$ and $CI$ (see Fig. \ref{fig2}(b)). If $M$ further increases to $M=N$, $CI$ agents invest in the first round but refuse to invest during the remaining rounds. In this case, each agent only retains the benefits of the first round of the game, and no one can obtain non-zero payoffs from the remaining rounds. The evolutionary outcomes in finite populations are similar to the result reported in Fig. \ref{fig2}(a), that is, the system spends most of the time in states where $T$ and $U$ agents coexist (see Fig. \ref{fig2}(c)).

Furthermore, we investigate how the investment threshold influences the stationary frequency of strategies $CI, T,$ and $U$, as shown in Fig. \ref{fig3}. We find that when the threshold $M$ is zero, $U$ can be more prevalent than $CI$ and $T$. With the increase of threshold, $CI$ and $T$ agents can form a strong alliance, preventing the invasion of $U$ agents. If $M$ further increases to $N$, $CI$ agents only invest in the first round, and thus all agents can only obtain benefit from the first round. In this case, the advantage of $U$ agents is greater than that of $T$ and $CI$ agents.

To probe deeper into the underlying mechanisms responsible for such a significant improvement of trust stemming from the investment threshold, we show in Fig. \ref{fig4} how the overall average payoffs of $CI, T$, and $U$ agents change with $M$. It can be seen that the average payoffs of all agents increase first and then decrease. Particularly, when $M=0$ or $M=N$, the average payoffs of $CI$ agents are negative (see the inset of Fig. \ref{fig4}). In general, an intermediate investment threshold can ensure a more evident advantage of alliance formed by $CI$ and $T$ agents over $U$ agents, even if this advantage can shrink with the increase of $M$.

\begin{figure*}[!h]
\begin{center}
\includegraphics[width=5.5in]{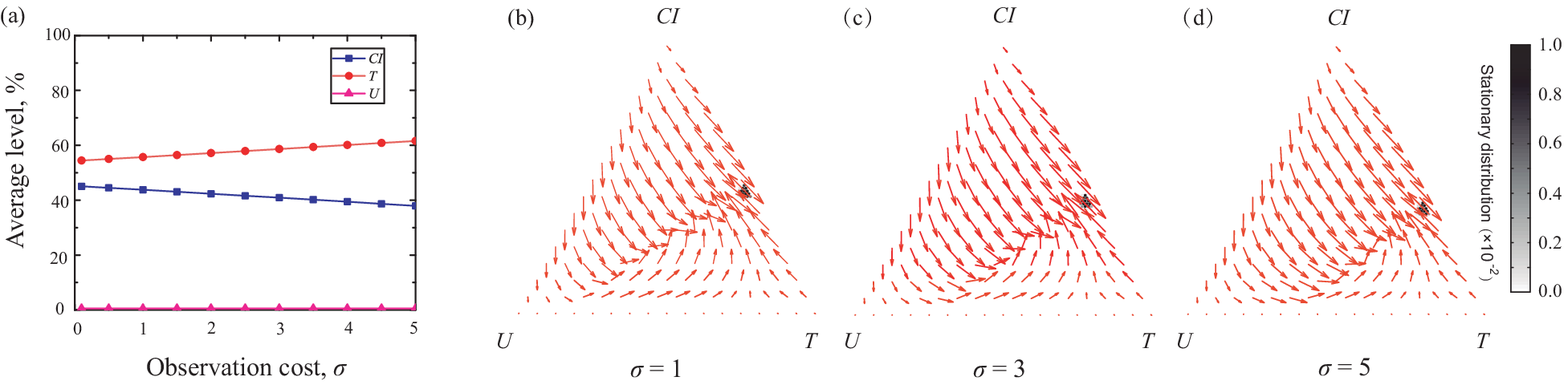}
\caption{Robustness of the evolutionary outcomes of $CI, T$, and $U$ to the changes of the observation cost $\sigma$. Panel (a) shows the average levels of three strategies as a function of $\sigma$. Panels (b)-(d) show the stationary distribution and the gradient of selection for three different $\sigma$ values. Parameter values are $Z=100, N=4, M=2, R_{U}=8, w=0.8, R_{T}=6, \mu=1/Z, \beta=5, tv=1$ in panel (a); $Z=100, N=4, M=2, tv=3, R_{U}=8, w=0.8, R_{T}=6, \mu=1/Z, \beta=5, \sigma=1$ in panel (b); $Z=100, N=4, M=2, tv=1, R_{U}=8, w=0.8, R_{T}=6, \mu=1/Z, \beta=5, \sigma=3$ in panel (c); $Z=100, N=4, M=2, tv=1, R_{U}=8, w=0.8, \sigma=5, R_{T}=6, \mu=1/Z, \beta=5$ in panel (d).}\label{fig7}
\end{center}
\end{figure*}

\begin{figure*}[!t]
\begin{center}
\includegraphics[width=5.5in]{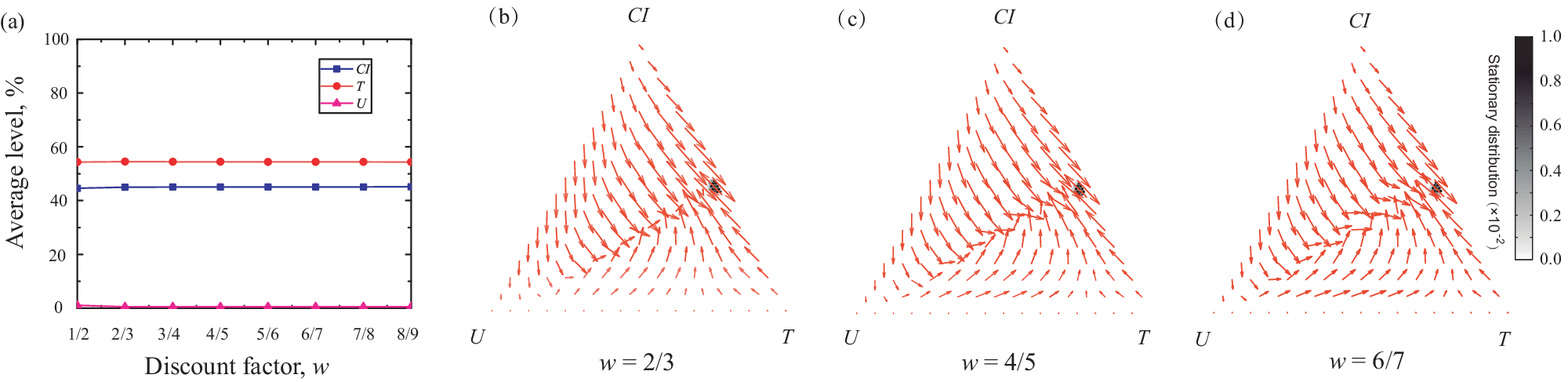}
\caption{Robustness of the evolutionary outcomes of $CI, T$, and $U$ to the changes of the discount factor $w$. Panel (a) shows the average levels of three strategies as a function of $w$. Panels (b)-(d) show the stationary distribution and the gradient of selection for three different $w$ values. Parameter values are $Z=100, N=4, M=2, R_{U}=8, \sigma=0.1, R_{T}=6, \mu=1/Z, \beta=5, tv=1$ in panel (a); $Z=100, N=4, M=2, tv=3, R_{U}=8, w=2/3, R_{T}=6, \mu=1/Z, \beta=5, \sigma=0.1$ in panel (b); $Z=100, N=4, M=2, tv=1, R_{U}=8, w=4/5, R_{T}=6, \mu=1/Z, \beta=5, \sigma=0.1$ in panel (c); $Z=100, N=4, M=2, tv=1, R_{U}=8, w=6/7, \sigma=0.1, R_{T}=6, \mu=1/Z, \beta=5$ in panel (d).}\label{fig8}
\end{center}
\end{figure*}

In order to illustrate the robustness of all results obtained in this paper, we first present how evolutionary dynamics of $CI, T$, and $U$ change with the multiply factor of $U$ agents, $R_{U}$ (see Fig. \ref{fig5}). Concretely, we show how the average levels of three strategies vary with $R_{U}$ in Fig. \ref{fig5}(a), and we find the average levels of $CI, T$, and $U$ remain almost constant with the increase of $R_{U}$. In addition, we show the evolutionary outcomes of the gradient of selection and the stationary distribution, which allow visualization of the dynamics, for three different values $R_{U}$ in Fig. \ref{fig5}(b)-(d). It shows that the population spends a significant time near the $CI$-$T$ edge of the simplex $S_{3}$, and most of the arrows in the simplex flow to the intermediate region of $CI$-$T$ edge, indicating that conditional investors can form an alliance with trustworthy trustees. It is worth pointing out that we set $R_{U}=6.66$, 7.98, and $9.96$, which can respectively correspond to the mild, moderate, and harsh interaction environment in the trust game (see Ref. \cite{chica2017networked}). Different from previous findings in \cite{chica2017networked}, our results show that trust can be maintained no matter whether the environment of the trust game is mild or harsh.

\begin{figure*}[!h]
\begin{center}
\includegraphics[width=5.5in]{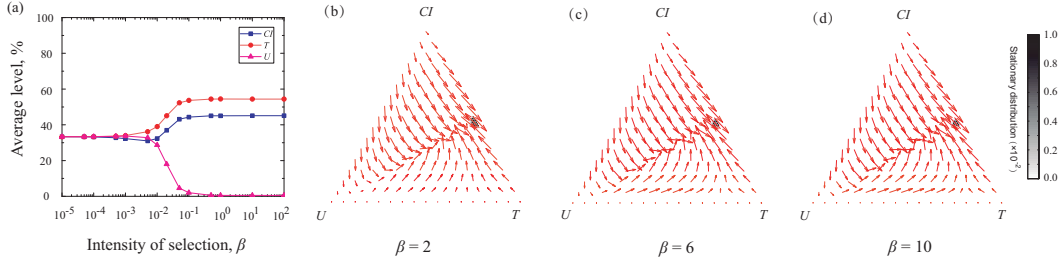}
\caption{Robustness of the evolutionary outcomes of $CI, T$, and $U$ to the changes of the intensity of selection $\beta$. Panel (a) shows the average levels of three strategies as a function of $\beta$. Panels (b)-(d) show the stationary distribution and the gradient of selection for three different $\beta$ values. Parameter values are $Z=100, N=4, M=2, R_{U}=8, \sigma=0.1, R_{T}=6, w=0.8, \mu=1/Z, tv=1$ in panel (a); $Z=100, N=4, M=2, tv=3, R_{U}=8, w=0.8, R_{T}=6, \mu=1/Z, \beta=2, \sigma=0.1$ in panel (b); $Z=100, N=4, M=2, tv=1, R_{U}=8, w=0.8, R_{T}=6, \mu=1/Z, \beta=6, \sigma=0.1$ in panel (c); $Z=100, N=4, M=2, tv=1, R_{U}=8, w=0.8, \sigma=0.1, R_{T}=6, \mu=1/Z, \beta=10$ in panel (d).}\label{fig9}
\end{center}
\end{figure*}

\begin{figure*}[!h]
\begin{center}
\includegraphics[width=5.5in]{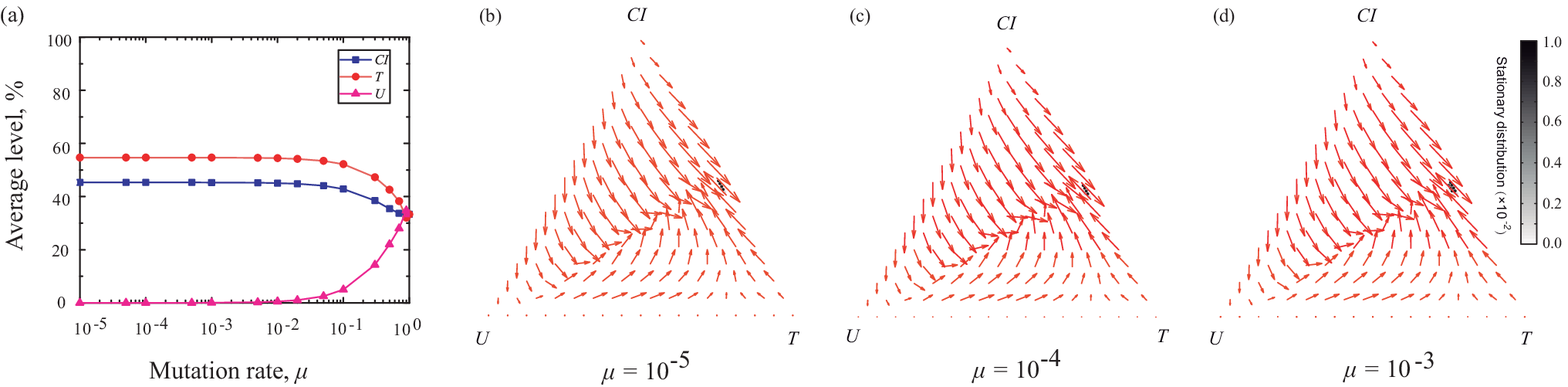}
\caption{Robustness of the evolutionary outcomes of $CI, T$, and $U$ to the changes of the mutation rate $\mu$. Panel (a) shows the average levels of three strategies as a function of $\mu$. Panels (b)-(d) show the stationary distribution and the gradient of selection for three different $\mu$ values. Parameter values are $Z=100, N=4, M=2, R_{U}=8, \sigma=0.1, R_{T}=6, w=0.8, \beta=5, tv=1$ in panel (a); $Z=100, N=4, M=2, tv=3, R_{U}=8, w=0.8, R_{T}=6, \mu=10^{-5}, \beta=5, \sigma=0.1$ in panel (b); $Z=100, N=4, M=2, tv=1, R_{U}=8, w=0.8, R_{T}=6, \mu=10^{-4}, \beta=5, \sigma=0.1$ in panel (c); $Z=100, N=4, M=2, tv=1, R_{U}=8, w=0.8, \sigma=0.1, R_{T}=6, \mu=10^{-3}, \beta=5$ in panel (d).}\label{fig10}
\end{center}
\end{figure*}

\begin{figure*}[!h]
\begin{center}
\includegraphics[width=5.5in]{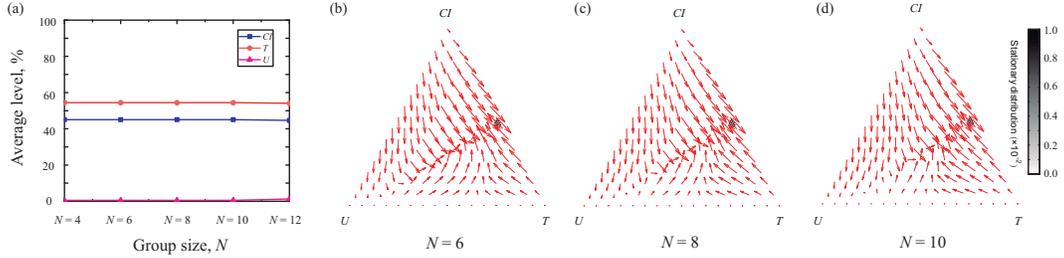}
\caption{Robustness of the evolutionary outcomes of $CI, T$, and $U$ to the changes of the group size $N$. Panel (a) shows the average levels of three strategies as a function of $N$. Panels (b)-(d) show the stationary distribution and the gradient of selection for three different $N$ values. Parameter values are $Z=100, M/N=0.5, R_{U}=8, \sigma=0.1, R_{T}=6, w=0.8, \beta=5, \mu=1/Z, tv=1$ in panel (a); $Z=100, N=6, M=3, tv=3, R_{U}=8, w=0.8, R_{T}=6, \mu=1/Z, \beta=5, \sigma=0.1$ in panel (b); $Z=100, N=8, M=4, tv=1, R_{U}=8, w=0.8, R_{T}=6, \mu=1/Z, \beta=5, \sigma=0.1$ in panel (c); $Z=100, N=10, M=5, tv=1, R_{U}=8, w=0.8, \sigma=0.1, R_{T}=6, \mu=1/Z, \beta=5$ in panel (d).}\label{fig11}
\end{center}
\end{figure*}

In what follows, we investigate the impact of the investment cost of investors on the evolutionary dynamics (see Fig. \ref{fig6}). We find that with the increase of $tv$ value, the average levels of three strategies are almost unchanged (see Fig. \ref{fig6}(a)). In addition, we find that changing the $tv$ values does not influence the gradient of selection and the stationary distribution of the system (see Fig. \ref{fig6}(b)-(d)).

In Fig. \ref{fig7}, we further present how observation cost affects the evolutionary dynamics of $CI$, $T$, and $U$ strategies. We find that the increase of observation cost will not change the evolutionary advantage of $CI$ and $T$ alliance over $U$ agents. Concretely, the average level of $CI$ decreases with the observation cost, the average level of $T$ increases, but the average level of $U$ remains almost unchanged (see Fig. \ref{fig7} (a)). Furthermore, we find that the entire population will spend most time near configurations in which $CI$ and $T$ agents coexist. In addition, with the increase of observation cost, the background shadow area is gradually moving towards the vertex $T$ (see Fig. \ref{fig7}(b)-(d)).

In Fig. \ref{fig8}, we investigate the role of discount factors $w$ in the evolutionary dynamics of $CI$, $T$, and $U$ strategies. Our results show that the average levels of these three strategies remain almost unchanged with the increase of $w$ (see Fig. \ref{fig8} (a)). Besides, for three different $w$ values, we find that our results are robust, that is, the population will spend most of the time in the intermediate
region of $CI-T$ edge, and most arrows flow to the coexistence states of $CI$ and $T$ (see Fig. \ref{fig8}(b)-(d)).

Furthermore, the effects of the mutation rate $\mu$ and the intensity of selection $\beta$ on the evolutionary dynamics are investigated (Figs. \ref{fig9} and \ref{fig10}).
When $\beta$ is small ($\beta<10^{-3}$) or when $\mu$ is significantly large (close to 1), the evolutionary process is mainly affected by imitation error or behavioral mutation,
which leads to the average levels of all strategies close to $1/3$ (see Figs. \ref{fig9}(a) and \ref{fig10}(a)). As $\beta$ increases or $\mu$ decreases, a strategy's
performance becomes increasingly important for the strategy's survival, eventually favoring the emerging alliance of the $CI$ and $T$ strategies, which crowds out the most $U$ agents (see Fig. \ref{fig9} (b)-(d) and Fig. \ref{fig10} (b)-(d)).

In what follows, we investigate how evolutionary dynamics of $CI, T$, and $U$ strategies change with the group size $N$ (Fig. \ref{fig11}). Here we set $M/N=0.5$, and we find that with the increase of $N$, the average levels of $CI, T$, and $U$ remain almost unchanged (see Fig. \ref{fig11}(a)). In particular, the average level of $T$ is the highest, the second largest frequency is formed by $CI$, while $U$ makes up the smallest fraction. In addition, the dynamic visualization results show that the population spends most time in the intermediate
region of $CI-T$ edge, and most arrows point to the states where $CI$ and $T$ agents coexist when $N$ changes appropriately (see Fig. \ref{fig11}(b)-(d)).

At the end of this section, we would like to point out that the work done here can have wider implications on other areas such as cognitive science, behavioral science, artificial intelligence, economics, and management science. Particularly, in the management system, developing and maintaining trust can promote the emergence of social exchange and economic transactions, which is significantly important for the effectiveness of management and organization. By investigating the evolutionary dynamics of trust in repeated group interactions, we show that the introduction of the conditional investment strategy can lead to the emergence of high trust and high trustworthiness in repeated group interactions, which can explain why trusting in others and reciprocating them with trustworthy actions can be sustained in social and economic interactions.

\section{Conclusions and Discussion}

In this work, we study the evolutionary dynamics of the conditional investment strategy in the NTG in which agents engage
in repeated group interactions. Different from the setting of NTG considered in previous works \cite{abbass2015n,hu2021adaptive}, we consider that investors can adjust their investment decisions according to the investment environment in repeated group interactions framework. Concretely, we use the tolerance threshold as an important parameter to describe the investment tendency of investors. We find that an intermediate threshold can lead to the evolution of high trust and high trustworthiness, while too low threshold level or too high threshold level cannot promote the evolution of trust. In summary, as the answers to the initially proposed questions, we can conclude that the introduction of the conditional investment strategy provides an avenue for trust to thrive in repeated group interactions.

As we have emphasized above, our model setup is well aligned with reality in which investors generally make decisions based on the current investment environment rather than investing blindly. A good example where our model could apply is the trust relationship between consumers and suppliers \cite{hawlitschek2016trust}. Consumers will refuse to consume if they find more suppliers selling inferior products during daily shopping on Amazon. Therefore, the collection of information is particularly important. The key assumption of conditional investment is that investors pay a permanent observation cost $\sigma$ to collect information about the trustworthiness of the trustee, and then make decisions in the subsequent possible group interactions based on what they observe. We have shown that increasing $\sigma$ leads to the decrease of the frequency of conditional investors and the increase of the frequency of trustworthy agents, while has little effect on the frequency of untrustworthy agents. Therefore, the result that the alliance formed by $CI$ and $T$ agents can effectively resist the invasion of $U$ agents is robust to the observation cost.

Lastly, it is worth emphasizing that in addition to the observation cost, our results are more robust to the change of other model parameters. Previous work on networked NTG has revealed that the level of trust is correlated with how ``difficult" the game is \cite{chica2017networked}. In our model, even though the social dilemma is extremely difficult (i.e. the ratio of temptation to defect $\frac{R_{U}}{R_{T}}>1.66$), a high level of trust can always be achieved (see Fig. \ref{fig5}(a)). Besides, for the sake of mathematical convenience, analysis of evolutionary dynamics of NTG has been mostly dealt with either in the limit of rare mutations \cite{lim2020stochastic} or in the limit of weak selection \cite{tarnita2015fairness}. Here our approach can be applied to arbitrary mutation rate and arbitrary intensity of selection values. We find that appropriate changes in mutation rate and intensity of selection will not affect our main results (see Fig. \ref{fig9} and Fig. \ref{fig10}).

Future work could explore the evolutionary dynamics of the repeated NTG in structured populations where the interactions among agents are typically not random but rather limited to a subset of the population \cite{chen2008pre,Szolnoki2015,zhang2021is}. Indeed, the network where agents interact with others has different structures, and the MDP approach adopted in our work is suitable for complete graph networks and can also be extended to other network structures \cite{chen10TC}. Besides, we can use theoretical approximation and agent-based simulations to study the evolutionary dynamics of trust on any population structure \cite{allen2017}.

\section*{CRediT authorship contribution statement}
Linjie Liu: Conceptualization, Methodology, Writing - original draft. Xiaojie Chen: Conceptualization, Formal
analysis, Writing - review \& editing, Supervision.

\section*{Declaration of Competing Interest}
The authors declare that they have no competing financial interests.

\section*{Acknowledgment}

This work was supported by the National Natural Science Foundation of China (Grants Nos. 61976048 and 62036002) and the Fundamental Research Funds of the Central Universities of China. L.L. acknowledges the support from Special Project of Scientific and Technological Innovation (Grant No. 2452022107).






\end{document}